\documentclass{amsart}

\usepackage[foot]{amsaddr}

\author{Matthew D. Kvalheim}
\address
{Department of Mathematics, University of Michigan,
	Ann Arbor, MI, USA}

\email{kvalheim@umich.edu}

\title[Obstructions to asymptotic stabilization]{Obstructions to asymptotic stabilization}

\usepackage{import} 

\usepackage[utf8]{inputenc}
\usepackage[T1]{fontenc}
\usepackage{lmodern}
\usepackage{amsmath}
\usepackage{amssymb}
\usepackage{amsthm}
\usepackage{stmaryrd} 
\usepackage{mathtools}
\usepackage{tikz-cd}
\usepackage{subfig}

\usepackage{thmtools}
\usepackage{thm-restate} 

\newcommand{\concept}[1]{\textbf{#1}}

\newcommand{\N}{\mathbb{N}}
\newcommand{\Z}{\mathbb{Z}}
\newcommand{\R}{\mathbb{R}}

\newcommand{\T}{T}

\newcommand{\id}{\textnormal{id}}

\newcommand{\Hom}{H}

\newcommand{\st}{M}
\newcommand{\att}{A}
\newcommand{\sa}{\st\ct\att}
\newcommand{\ua}{U\ct\att}

\newcommand{\vp}{\varphi}

\newcommand{\ct}{\setminus}
\newcommand{\tX}{\tilde{X}}
\newcommand{\tY}{\tilde{Y}}
\newcommand{\tpX}{\Phi_{\tilde{X}}}
\newcommand{\tpY}{\Phi_{\tilde{Y}}}
\newcommand{\Exp}{\textnormal{Exp}}

\theoremstyle{definition}

\newtheorem*{Lem-non}{Lemma}
\newtheorem{Th}{Theorem}

\newtheorem*{Co-non}{Corollary}

\newtheorem*{Prop-non}{Proposition}

\newtheorem*{Quest-non}{Question}

\newtheorem*{Th-non}{Theorem}

\newcommand{\thistheoremname}{}
\newtheorem*{genericthm}{\thistheoremname}
{\renewcommand{\thistheoremname}{Theorem~\ref{#1}$'$}%
	\begin{genericthm}}
	{\end{genericthm}}

\newtheorem*{Def*}{Definition}

\newtheorem*{Ex-non}{Example}

\newtheorem{Rem}{Remark}
\newtheorem*{Note}{Note}

\usepackage{marvosym}
\usepackage[
hypertexnames,%
citecolor=blue,%
colorlinks=true,%
linkcolor=red%
]{hyperref}
\usepackage[all]{hypcap}

\begin{document}
	
\begin{abstract}	
Necessary conditions for asymptotic stability and stabilizability of  subsets for dynamical and control systems are obtained.
The main necessary condition is homotopical and is in turn used to obtain a homological one.
A certain extension is ruled out. 
Questions are posed.
\end{abstract}

\maketitle

\section{Intrigue}\label{sec:intrigue}
The following is a simplified version of the main result, Theorem~\ref{th:dist}.

\begin{Th-non}
Let $X$, $Y$ be smooth vector fields on a manifold $\st$ with a compact set $\att\subset \st$ asymptotically stable for both.
There is an open set $U\supset \att$ such that $X|_{U\ct\att}$, $Y|_{U\ct\att}$ are homotopic through nowhere-zero vector fields.
\end{Th-non}
Motivation is given in sections~\ref{sec:introduction} and \ref{sec:results}.
Briefly, this result for \emph{dynamical} systems yields feedback stabilizability ``tests'' for \emph{control} systems generalizing and/or strengthening some in the literature.
Theorem~\ref{th:dist} enables further such tests via Theorem~\ref{th:control-homology} or obstruction theory \cite[Part~III]{steenrod1999topology}.
Open questions  are in section~\ref{sec:questions}.


\section{Introduction}\label{sec:introduction}
In control theory one studies ordinary differential equations 
\begin{equation}\label{eq:cs-intro}
\frac{dx}{dt} = f(x,u) 
\end{equation}
depending on a control parameter $u$ where $x$ lives in a smooth manifold $\st$ and $f(x,u)\in \T_x \st$ for all $x$ and $u$.
Given a compact subset $\att$ of $\st$, the \concept{feedback stabilization problem} concerns proving (non)existence and ideally constructing a \concept{feedback law} $x\mapsto u(x)$ of suitable regularity such that $\att$ is asymptotically stable for the \concept{closed-loop vector field} $x\mapsto f(x,u(x))$. 
This paper merely considers proving nonexistence.

In the special case that $\att$ is a point and $\st$ is a Euclidean space, \concept{Brockett's necessary condition} says in particular that the existence of a $C^1$ closed-loop vector field solving the feedback stabilization problem implies that the image of the map $(x,u)\mapsto f(x,u) \in \R^n$ contains a neighborhood of $0$ \cite[Thm~1.(iii)]{brockett1983asymptotic}. 
(Here $f$ is viewed as $\R^n$-valued via the usual tangent bundle identification $\T \R^n \approx \R^n \times \R^n$.) 
To achieve this, Brockett derived a necessary condition satisfied by a stabilizing closed-loop vector field---a \emph{dynamical} system---and parlayed this into one for \emph{control} systems \eqref{eq:cs-intro}.
This technique is used in later work and in this paper, where the main result for dynamics (Theorem~\ref{th:dist}) implies one for control (Theorem~\ref{th:control-homology}).

The monograph by Krasnosel'ski\u{\i} and Zabre\u{\i}ko appeared around the same time as Brockett's paper and has been influential in this area. 
It contains necessary conditions for asymptotic stability of an equilibrium point of a vector field involving its index \cite[sec.~52]{krasnoselskii1984geometrical}.
These results were used in part by Zabczyk to give an alternative proof and extension of Brockett's necessary condition \cite{zabczyk1989some}.
Further necessary conditions for feedback stabilizability have since been discovered for the case that $\att$ is a point \cite{coron1990necessary,christopherson2022feedback} or a closed submanifold of a Euclidean space \cite{mansouri2007local,mansouri2010topological} or a compact subset of a smooth manifold \cite{kvalheim2022necessary}.\footnote{Exponential \cite{gupta2018linear}, global \cite{byrnes2008brockett,baryshnikov2021topological}, time-varying 
\cite{coron1992global}, and discontinuous  
\cite{clarke1997asymptotic} variants of the feedback stabilization problem are not considered in the present paper.}

Such necessary conditions yield practical ``tests''  relieving the expenditure of time and resources searching for stabilizing feedback laws that do not exist.
These tests are easy to apply in a variety of concrete instances of \eqref{eq:cs-intro}, and the author hopes that the homological tests among these (\cite{coron1990necessary,mansouri2007local,mansouri2010topological}, Theorem~\ref{th:control-homology}) might be usefully numerically automatable \cite{kaczynski2004computational}.
The purpose of the present paper is to introduce necessary conditions providing more general and/or stronger tests.

This is done in section~\ref{sec:results}.
The main result is proved in section~\ref{sec:proof}.
Some questions are posed in section~\ref{sec:questions}.

\section{Results}\label{sec:results}
A $C^0$ vector field is \concept{uniquely integrable} if it has unique maximal integral curves.\footnote{Such vector fields generate $C^0$ local flows (see \cite[App.~A.1]{kvalheim2022necessary} for details).
Standard examples are $C^\infty$ or $C^1$ or even locally Lipschitz vector fields.}
A compact invariant set $\att$ for such a vector field is \concept{asymptotically stable} if for every open set $W\supset \att$ there is an open set $V\supset \att$ such that all maximal forward integral curves initialized in $V$ are defined for all positive time, contained in $W$, and converge to $\att$.
\begin{restatable}[]{Th}{ThmDist}\label{th:dist}
Let $X$, $Y$ be uniquely integrable $C^0$ vector fields on a $C^\infty$ manifold $\st$ with a compact set $\att\subset \st$ asymptotically stable for both.
There is an open set $U\supset \att$ such that $X|_{\ua}$, $Y|_{\ua}$ are homotopic through nowhere-zero vector fields.
\end{restatable}

\begin{Note}
$U$ admits an explicit description  (Remark~\ref{rem:u-explicit}).
\end{Note}

When $\att$ is a point, Theorem~\ref{th:dist} specializes to a known result  \cite[p.~340]{krasnoselskii1984geometrical} which implies Brockett's necessary condition \cite{zabczyk1989some}.

When $\att$ is known to be asymptotically stable for $Y$, Theorem~\ref{th:dist} can be used to test whether $\att$ is  also so for $X$, e.g. by using obstruction theory \cite[Part~III]{steenrod1999topology}.
If $\sa$ is orientable and $(n+1)$-dimensional, the obstructions to $X|_{\ua}$, $Y|_{\ua}$ being homotopic through nowhere-zero vector fields lie in cohomology groups with coefficients in the homotopy groups of the $n$-sphere \cite[Thm~1.1]{davis1972vector}.

Alternatively, Theorem~\ref{th:dist} directly implies the following homological result that can also be used to test if $\att$ is asymptotically stable for $X$ assuming it is so for $Y$.
This result will be used to formulate a necessary condition for feedback stabilizability of control systems in Theorem~\ref{th:control-homology}.

\begin{Lem-non}
Let $X$, $Y$ be uniquely integrable $C^0$ vector fields on a $C^\infty$ manifold $\st$ with a compact set $\att\subset \st$ asymptotically stable for both.
For all sufficiently small open sets $U\supset \att$ the induced graded homomorphisms
$$(X|_{\ua})_*, (Y|_{\ua})_* \colon  \Hom_{\bullet}(\ua)\to \Hom_{\bullet}(\T(\ua)\ct 0)$$
on singular homology coincide.
\end{Lem-non}
Control systems will now now be defined more formally than in section~\ref{sec:introduction}.
The definitions here are somewhat less general than in \cite[sec.~3]{kvalheim2022necessary} (continuity assumptions are added on $E$, $p$, $u$).
Brockett's definition \cite{brockett1977control} motivates the following one used in this paper: a \concept{control system} is a tuple $(E, \st, p, f)$ of a topological space $E$, a $C^\infty$ manifold $\st$, a continuous surjection $p\colon E\to \st$, and a fiber-preserving map $f\colon E\to \T \st$ (i.e. $q \circ f = p$ where $q\colon \T \st\to \st$ is the tangent bundle projection).
The standard example in control theory is $E = \st\times \R^m$ for some $m\in \N$ with $p(x,u) = x$.
In this paper a \concept{feedback law} $u\colon \st \to E$ is a continuous section of $p$ (i.e. $p\circ u = \id_\st$) such that the $C^0$ vector field $f\circ u$ is uniquely integrable. Note: $f\circ u$ is indeed a vector field on $\st$ since $q\circ f\circ u = p\circ u = \id_\st$.

The lemma implies the next result by taking $X = f\circ u$ and using $X_* = f_* \circ u_*$.

\begin{Th}\label{th:control-homology}
Let $\st$ be a $C^\infty$ manifold, $(E,\st,p,f)$ be a control system, and $\att\subset \st$ be a compact subset that is asymptotically stable for \emph{some} uniquely integrable $C^0$ vector field $Y$ on $\st$.
If there is a feedback law $u$ such that $\att$ is asymptotically stable for $f\circ u$, then
for all sufficiently small open sets $U \supset \att$,
\begin{equation}\label{eq:control-homology}
Y_*\Hom_{\bullet}(\ua) \subset f_*\Hom_{\bullet}(p^{-1}(\ua)\ct f^{-1}(0)) \subset \Hom_{\bullet}(\T(\ua) \ct 0),
\end{equation}
where $Y_*$, $f_*$ are the induced graded homomorphisms on singular homology.
\end{Th}

\begin{Rem}\label{rem:coron-mansouri}
Theorem~\ref{th:control-homology} generalizes results of Coron \cite[Thm~2]{coron1990necessary} and Mansouri \cite[Thm~4]{mansouri2007local}, \cite[Thm~2.3]{mansouri2010topological} in a few ways.
E.g. Coron assumes $\att$ is a point, Mansouri assumes $\att$ is a compact boundaryless submanifold, and both assume $\st = \R^n$.
These assumptions are not made here.
However, the conclusions of Coron's and Mansouri's theorems are more explicit.
\end{Rem}

The following example illustrates Theorems~\ref{th:dist}, \ref{th:control-homology}.
Its conclusion is not deducible from the results of \cite{brockett1983asymptotic,zabczyk1989some,coron1990necessary} since $\att$ is not a point, from \cite{kvalheim2022necessary} since $\att$ has zero Euler characteristic, or from \cite{mansouri2007local,mansouri2010topological} for the same reason or because $\st$ does not embed as an open subset of a Euclidean space. 
\begin{Ex-non}
	Let $\st$ be the M\"{o}bius band viewed as the quotient of $\R^2$ by the identifications $(x,y)\sim (x+2\pi,-y)$.
	For $\varepsilon > 0$ the sets $\{y=0\}$, $\{|y|=\varepsilon\}$ and vector fields $\partial_x$,  $-y\partial_y$ descend to well-defined sets $\att, C_\varepsilon \subset \st$ and vector fields $\partial_\theta$, $Y$ on $\st$, where $\theta = x\mod 2\pi$.
	Is there a $C^0$ feedback law $u\colon \st \to \R$ such that $\att$ is asymptotically stable for the uniquely integrable $X \coloneqq (\cos \theta \partial_\theta + \sin \theta Y)u$?
	
	If so, fix an open set $U\supset \att$ small enough that $X|_{\ua}$ is nowhere-zero and $\varepsilon > 0$ small enough that $C_\varepsilon \subset U$.
	The vector fields $X|_{\ua}$, $Y|_{\ua}$ are not homotopic since, with respect to the frame $(\partial_\theta, Y)|_{C_\varepsilon}$ for $\T \st|_{C_\varepsilon}$,  $X|_{C_\varepsilon}$ has winding number $2$  while $Y|_{C_\varepsilon}$ has winding number $0$.
	Since $\att$ is asymptotically stable for $Y$ this contradicts Theorem~\ref{th:dist}, so the desired $u$ does not exist. 
	
	This conclusion is also detected by Theorem~\ref{th:control-homology}.
	For this example the control system is $(\st \times \R, \st, p, f)$ where $p(m,u)=m$ and $f(m,u) = (\cos\theta \partial_\theta + \sin \theta Y)u$.
	For each $\varepsilon > 0$ the set $\{|y|<2\varepsilon\}$ descends to an open set $U\subset \st$ such that $\ua$ deformation retracts onto $C_\varepsilon$, and 
	the frame $(\partial_\theta, Y)|_{\ua}$ induces a diffeomorphism $\T(\ua) \ct 0 \approx (\ua)\times (\R^2\ct 0)$.
	Hence the first homology of $\T(\ua) \ct 0 $ is isomorphic to $\Z\oplus \Z$  with, by the observations above, generators $(e_1,e_2)$ such that $$Y_* H_1(\ua) = \Z e_1 \quad \textnormal{and} \quad f_* H_1(p^{-1}(\ua)\ct f^{-1}(0)) = \Z(e_1+2e_2).$$
	Since $\Z e_1 \ni e_1 \not\in \Z(e_1+2e_2)$, Theorem~\ref{th:control-homology} also implies that $u$ does not exist.
\end{Ex-non}

Theorems~\ref{th:dist}, \ref{th:control-homology} are not always very informative in spite of Remark~\ref{rem:coron-mansouri}.
For example, if $\att$ is the image of a periodic orbit with the same orientation for smooth vector fields $X$ and $Y$,  the straight-line homotopy over a sufficiently small open set  $U\supset \att$ satisfies the conclusion of Theorem~\ref{th:dist}  regardless of whether $\att$ is attracting or repelling or neither for $X$ or $Y$.  
Since Theorem~\ref{th:dist} is at least as strong as Theorem~\ref{th:control-homology}, neither result provides any information on stability in this case.
This is related to the fact that the Euler characteristic of a circle is zero (cf. \cite{mansouri2007local,mansouri2010topological,kvalheim2022necessary}).

This motivates a search for conclusions stronger than that of Theorem~\ref{th:dist} given vector fields $X$, $Y$ satisfying its hypotheses.
For example, can the homotopy of Theorem~\ref{th:dist} be taken through vector fields stabilizing $\att$?
This question is related to one of Conley and not answered in this paper but see section~\ref{sec:questions} for additional thoughts.
Alternatively, if $X$, $Y$ stabilize $\att$ and also have additional structure, one can ask if the homotopy of Theorem~\ref{th:dist} can be taken through vector fields preserving this structure.

For example, if $X$, $Y$ are gradient vector fields stabilizing $\att$, can the homotopy of Theorem~\ref{th:dist} be taken through gradient vector fields?
Unfortunately, the answer is ``yes'' but does not yield a stronger conclusion since the following proposition says that gradients on an open set are always homotopic through nowhere-zero gradients if they are homotopic through nowhere-zero vector fields.

\begin{Prop-non}
Let $W$ be an open subset of a $C^\infty$ Riemannian manifold $\st$, and let $F, G\colon W\to \R$ be $C^1$ functions such that $\nabla F$, $\nabla G$ are homotopic through nowhere-zero vector fields. 
Then there is also a $C^1$ homotopy from $F$ to $G$ through functions whose gradients are nowhere-zero.
\end{Prop-non}

\begin{proof}
Smooth approximation techniques \cite{hirsch1976differential} imply that the submersions $F$, $G$ are homotopic through $C^1$ submersions to $C^\infty$ submersions.
By taking metric duals the assumptions (and smooth approximation techniques again) imply that the differentials of these $C^\infty$ submersions are $C^\infty$ homotopic through nowhere-zero differential one-forms.
Since $W$ is open and submersions satisfy Gromov's $h$-principle for open $\textnormal{Diff}(W)$-invariant differential relations \cite{gromov1969stable} (lucidly explained in \cite[Thm~7.2.3]{eliashberg2002intro}), the $C^\infty$ submersions are $C^\infty$ homotopic through submersions (which have nowhere-zero gradients by definition).	
\end{proof}

\begin{Note}
	When $W = \ua$, the homotopy of the proposition need not be through functions whose gradients ``point toward'' $\att$  (cf. \cite[Fig.~4.1]{eliashberg2002intro} and the question about homotopies through stabilizing vector fields above and in section~\ref{sec:questions}).
\end{Note}

\section{Proof of Theorem~\ref{th:dist}}\label{sec:proof}
The proof of Theorem~\ref{th:dist} uses the following.
If $Y$ is a uniquely integrable $C^0$ vector field with an asymptotically stable compact set $\att$ having basin of attraction $B$, a \concept{proper $C^\infty$ Lyapunov function} is a $C^\infty$ function $V\colon B\to [0,\infty)$ with $\att = V^{-1}(0)$, $\{V \leq c\}$ compact for all $c < \infty$, and Lie derivative $YV < 0$ on $B\ct \att$.

\ThmDist*
\begin{Rem}\label{rem:u-explicit}
The proof yields an explicit description of one such $U$.
Let $B_X, B_Y \subset \st$ denote the basins of attraction of $\att$ with respect to $X$, $Y$.
Then $U$ can be taken to be the largest subset of $B_X\cap B_Y$ such that all forward $X$-orbits with initial condition in $U$ are contained in $B_Y$.
Equivalently, $U$ is the basin of attraction of $\att$ for $X|_{B_Y}$, so $U$ is open in $B_Y$ and hence also in $\st$.
\end{Rem}

\begin{proof}
Without loss of generality $\st$ is redefined to be the basin of attraction of $\att$ for $Y$.
Let $U\subset \st$ be the basin of attraction of $\att$ for $X$ and $\vp\colon \st \to (0,\infty)$ be a $C^\infty$ function such that the vector fields $\vp X$, $\vp Y$ are complete.
Since $\att$ is asymptotically stable for these complete vector fields there exist proper $C^\infty$ Lyapunov functions $V_X\colon U\to [0,\infty)$, $V_Y\colon \st \to [0,\infty)$ for $\vp X$, $\vp Y$ \cite[Thm~3.2]{wilson1969smooth}, \cite[sec.~6]{fathi2019smoothing}.
Fix any $C^\infty$ Riemannian metric and let $\psi\colon \st\to (0,\infty)$ be a $C^\infty$ function such that $- \psi \nabla V_X$, $-\psi \nabla V_Y$ are complete.
Note that $\att$ is asymptotically stable for the latter pair of vector fields with basins of attraction $U$, $\st$, and the straight-line homotopies from $X$ to $-\psi \nabla V_X$ and $Y$ to $-\psi \nabla V_Y$ are nowhere-zero on $\ua$ since $XV_X, Y V_Y < 0$ on $\ua$.
Thus, by replacing $X$, $Y$ with $- \psi \nabla V_X$, $-\psi \nabla V_Y$ it may be assumed from here on that $X$, $Y$ are complete and $C^\infty$ with flows $\Phi_X, \Phi_Y\in C^\infty$.

Asymptotic stability implies that, for each $z\in \ua$, there is $t_z > 0$ and an open subset $U_z\ni z$ of $\ua$ such that $\Phi_{X}^{[t_z,\infty)}(m)\cap \Phi_Y^{(-\infty,0]}(m)=\varnothing$ for all $m\in U_z$.\footnote{Proof: let $c = V_Y(z)$. 
Asymptotic stability of $\att$ for $X$ implies the existence of an open set $W$ and $t_z > 0$ such that $\att \subset W\subset \{V_Y < c/2\}$, $\Phi_X^{t_z}(z)\in W$, and $\Phi_X^{[0,\infty)}(W)\subset \{V_Y < c/2\}$.
By continuity there is an open set $U_z\ni z$ such that $U_z\subset \{V_Y > 2c/3\}$ and $\Phi_X^{t_z}(U_z)\subset W$, so $\Phi_X^{[t_z,\infty)}(U_z)\subset \Phi_X^{[0,\infty)}(W) \subset \{V_Y<c/2\}$ and $\Phi_Y^{(-\infty,0]}(U_z)\subset \{V_Y > 2c/3\}$ are disjoint.}
Let $(\psi_z)_{z\in \ua}$ be a $C^\infty$ partition of unity \cite[Thm~2.23]{lee2013smooth} subordinate to the open cover $(U_z)_{z\in \ua}$ of $\ua$ and define the $C^\infty$ function $\tau\colon \ua\to (0,\infty)$ by $\tau(m)\coloneqq \sum_{z\in \ua} t_z \psi_z(m)$.
Then for each $m\in \ua$ there is $z\in \ua$ such that $m\in U_z$ and $\tau(m)\geq t_z$.
Hence $\tau$, to be used later, satisfies
\begin{equation}\label{eq:tau-property}
\textnormal{$\Phi_X^{[\tau(m),\infty)}(m)\cap \Phi_Y^{(-\infty,0]}(m) = \varnothing$ for all $m\in \ua$.}
\end{equation}

Next, fix any $c\in (0,\infty)$ and define $N\coloneqq V_Y^{-1}(c)$.
Then 
(cf. \cite[Thm~3.2]{wilson1967structure})
\begin{equation}\label{eq:diffeo-to-product}
\textnormal{$\Phi_Y|_{\R\times N}\colon \R\times N \approx \sa$ is a diffeomorphism pushing forward $1\times 0$ to $Y$},
\end{equation}
so it may be assumed from here on that $\sa = \R\times N$ and $Y=1\times 0$.

Equip $N$ with any $C^\infty$ Riemannian metric, let $\exp^N\colon \T N\to N$ be the corresponding exponential map, and let $\exp^\R\colon \T\R\to \R$ be the exponential map of the standard Euclidean metric.
Define $$\exp\coloneqq \exp^\R \times \exp^N\colon \T(\sa)\to \sa$$ and let $\pi\colon \T (\sa)\to \sa$,  $\pi_\R\colon \T \R\to \R$, $\pi_N\colon \T N\to N$ be the tangent bundle projections.
By a standard argument $(\pi_N,\exp^{N})\colon \T N\to N^2$ maps a neighborhood of the zero section $0_{\T N}$ diffeomorphically onto its image (see \cite[p.~121]{brocker1982intro}), and clearly $(\pi_{\R}, \exp^{\R})\colon \T \R\to \R^2$ is a diffeomorphism. 

It follows in particular that $\Exp\coloneqq (\pi,\exp)|_{\T (\ua)}\colon \T(\ua)\to (\ua)\times (\sa)$ is a $C^\infty$ diffeomorphism from a neighborhood $U_0$ of $(\T \R\times 0_{\T N})|_{\ua}$ onto its open image $U_1$.
Moreover, by construction $\exp(tY(m)) = \exp_m(t\times 0) =  \Phi_Y^t(m)$ for every $m\in \ua$.
Defining $\tpX$, $\tpY$ to be the flows of $\tX\coloneqq (0\times X)$, $\tY\coloneqq (0\times Y)$ on $(\ua)\times (\sa)$ and $\Delta\coloneqq \{(m,m)\colon m\in \ua\}$, it follows that  $\Phi_{\tY}^{(-\infty,0]}(\Delta)\subset U_1$ is a closed subset, which by \eqref{eq:tau-property} is disjoint from $ \bigcup_{m\in \ua}\tpX^{[\tau(m),\infty)}(m,m)$.
Moreover, the corresponding version of \eqref{eq:diffeo-to-product} for $X$ rather than $Y$ implies that the latter subset is closed in $(\ua)\times (\sa)$.

Thus, there exists a $C^\infty$ function $\alpha\colon (\ua)\times (\sa)  \to [0,1]$ equal to $0$ on neighborhoods of $ \bigcup_{m\in \ua}\tpX^{[\tau(m),\infty)}(m,m)$ and $[(\ua)\times(\sa)]\ct U_1$, and equal to $1$ on a neighborhood of 
$ \tpY^{(-\infty,0]}(\Delta)$.
Define a $C^\infty$ function $\beta\colon [0,1]\times (\ua)\to [0,1]$ by $\beta(t,m)\coloneqq \alpha(\tpX^{t\tau(m)}(m,m))$.
Define $H\colon [0,1]\times (\ua)\to \T(\ua)\ct 0$ by
\begin{equation*}
\begin{split}
H(t,m)
&\coloneqq \begin{cases}
X(m), & t=0\\
Y(m), &  \beta(t,m) = 0\\
\frac{\beta(t,m)}{t\tau(m)}\Exp^{-1}(\tpX^{t\tau(m)}(m,m)) + (1-\beta(t,m))Y(m), & \textnormal{otherwise}.
\end{cases}
\end{split}
\end{equation*}
Clearly $H$ is continuous on $(0,1]\times (\ua)$, and continuity of $H$ at points in $\{0\}\times (\ua)$ follows from an application of Taylor's theorem in local coordinates.
Thus, $H$ is a homotopy from $X|_{\ua}$ to $Y|_{\ua}$ through nowhere-zero vector fields. 
\end{proof}

\section{Questions}\label{sec:questions}
In this section some remaining questions are posed and discussed.

Can the assumption of unique integrability be removed, as in \cite{orsi2003necessary}?

Motivated by the discussion before the proposition of section~\ref{sec:results}, can anything analogous be said for other special classes of vector fields that might be useful for strengthening the conclusion of Theorem~\ref{th:dist} under more specific hypotheses?

It was asked after the example of section~\ref{sec:results} whether the homotopy of Theorem~\ref{th:dist} can be taken through vector fields rendering $\att$ asymptotically stable.
It seems to the author that, using the $s$-cobordism theorem due to Barden-Mazur-Stallings \cite{barden1963structure,mazur1963relative,stallings1965on}, some sufficient conditions for this can be given when $\dim \st \geq 6$ under certain hypotheses on fundamental groups related to Whitehead torsion \cite{milnor1966whitehead}.
However, this would only provide a partial answer to the question.

It will now be explained how the question just mentioned relates to one asked by Conley.
In his monograph Conley defined what is now called the Conley index and showed that isolated invariant sets related by continuation (defined in \cite{conley1978isolated}) have ``the same'' Conley index.
Conley asked to what extent the converse is true \cite[Sec.~IV.8.1.A]{conley1978isolated}.
Simple examples show that isolated invariant sets with the same Conley index need not be related by continuation \cite[pp.~1--2]{mrozek2000conley}, as Conley surely knew, but one can still ask whether a converse holds in certain situations.
One such situation is when the pair of isolated invariant sets are the same asymptotically stable set but for two different vector fields; a homotopy of these through stabilizing vector fields would induce a continuation, and this explains the relationship of Conley's question with that of the previous paragraph.
Reineck has already investigated Conley's question using techniques related to the $h$-cobordism theorem (special case of the $s$-cobordism theorem) \cite{reineck1992continuation}. 

Finally, from the view of the present paper, it seems to the author that the ultimate question is the following.
Is there a ``universal'' necessary \emph{and} sufficient condition for stabilizability that can be used to actually \emph{construct} smooth stabilizing feedbacks via numerical means?
It seems to the author that a negative answer to this question would also be interesting.

    \section*{Acknowledgments}
    This work was funded by ONR N00014-16-1-2817, a Vannevar Bush Faculty Fellowship  sponsored by the Basic Research Office
    of the Assistant Secretary of Defense for Research and
    Engineering.
    The author gratefully acknowledges  Yu. Baryshnikov, B. A. Christopherson, J. Gould, R. Gupta, D. E. Koditschek, B. S. Mordukhovich, and T. O. Rot for helpful conversations and insights and references.

   	\bibliographystyle{amsalpha}
   	\bibliography{ref}
   	
\end{document}